\documentclass[12pt]{article}
\usepackage{amsmath,amssymb,amsthm,amsfonts}
\usepackage{latexsym}
\begin{document}

\newtheorem{prop}{Proposition}[section]
\newtheorem{cor}{Corollary}[section] 
\newtheorem{theo}{Theorem}[section]
\newtheorem{lem}{Lemma}[section]
\newtheorem{rem}{Remark}[section]

\renewcommand{\theequation}{\thesection.\arabic{equation}} 
\setcounter{page}{1} 
\noindent
\begin{center}
{\Large \bf Merging percolation and classical random graphs:
  Phase transition in dimension 1}
\end{center}

\begin{center}
TATYANA S. TUROVA\footnote{Research was 
supported by the Swedish Natural Science Research
Council.} and THOMAS VALLIER
\end{center}

\begin{center}
{\it Mathematical Center, University of
Lund, Box 118, Lund S-221 00, 
Sweden. }
\end{center}

\begin{abstract}
We study a  random graph model which combines properties of
the edge percolation model on $Z^d$ and a classical random graph $G(n,
c/n)$. We show that
this model, being a {\it homogeneous} random
graph,
has a natural relation to the so-called
"rank 1 case" of {\it inhomogeneous} random
graphs. This allows us to use the newly developed theory of
inhomogeneous random graphs to describe completely the phase diagram in
the case $d=1$.
The phase transition is similar to the classical random graph, it is
of the second order. We also find the scaled size of the largest connected
component above the phase transition.

\end{abstract}

\section{Introduction.}
\setcounter{equation}{0}

We consider a graph on the set of vertices $V_N^d:=\{1, \ldots, N\}^d$
in $Z^d$,
where the edges between any two different vertices $i$ and $j$ are
presented independently with probabilities
\begin{displaymath}
    p_{ij} = \left\{ \begin{array}{ll}
            q, & \mbox{ if } |i-j|= 1, \\
            c/N^d, & \mbox{ if } |i-j|>1 ,
        \end{array} \right.
\end{displaymath}
where $0\leq q \leq 1$ and $0<c<N$ are constants.
This graph, call it $G_N^d(q,c)$ is a mixture of percolation model,
where
each pair of
neighbours in $Z^d$ is connected with probability $q$, and a random
graph model, where each vertex is connected to any other vertex with
probability $c/|V_N^d|$.

The introduced model is a simplification of the most
common graphs designed to study natural  phenomena, in particular,  biological 
neural networks \cite{TV}.
Observe the difference between this and the so-called
"small-world" models intensively studied after \cite{WS}. In the 
 "small-world" models where edges from the grid  may be kept or
 removed,  only finite number (often at most  $2d$) of the long-range edges may
 come out of each vertex, and the probability of those is a  fixed
 number.
 
We are interested in the limiting behaviour of the introduced graph $G_N^d(q,c)$ as $N
\rightarrow \infty$. One can consider this model as a graph on $Z^d$ or on
a torus, in the limit the results are the same.
The one-dimensional case which we study here, is exactly
solvable. We shall write $G_N^1(q,c)=G_N(q,c)$.

Let $X$ be a random variable with
$Fs(1-q)$-distribution, i.e.,
\begin{equation}\label{X}
{\bf P} \left\{ X=k \right\} =(1-q)q^{k-1}, \ \ \ \ k=1,2,
  \ldots \ ,
\end{equation}
with
\[{\bf E} X = \frac{1}{1-q} \ .\]
Let further $ C_1 \Big( G \Big)$ denote the size of the largest
connected component in a graph $G$. 

\begin{theo}\label{T}
For any  $0\leq q < 1$ 
define
\begin{equation}\label{cr1}
  c^{cr} (q) = \frac{{\bf E} X }{{\bf E} X^2 } = \frac{1-q}{1+q} \ .
\end{equation}

\noindent
$\mbox{i)}$ If 
$c  < c^{cr} (q)$ then there exists a constant
$\alpha = \alpha (q,c)>\frac{1}{|\log q|}$ 
such that
\begin{equation}\label{E1}
 {\bf P} \left\{ C_1  \Big( G_N(q,c)  \Big) >  \alpha\, \log N 
\right\}  \rightarrow  0 ,
\end{equation}
and for any $\alpha _1<\frac{1}{|\log q|}$
\begin{equation}\label{E2}
{\bf P} \left\{ C_1  \Big( G_N(q,c)  \Big) < \alpha _1  \log N
\right\}  \rightarrow  0,
\end{equation}
 as  $N
\rightarrow \infty$.

\medskip

\noindent
$\mbox{ii)}$ If 
$c  \geq  c^{cr} (q)$ then 
\begin{equation}\label{lp}
\frac{C_1 \Big(
G_N(q,c)
    \Big)}{N}
\stackrel{P}{\rightarrow} \beta  
\end{equation}
as $ N \rightarrow \infty$,
with $\beta = \beta (q,c)$
defined as the maximal solution to
\begin{equation}\label{be}
  \beta = 1-\frac{1}{{\bf E} X}\ {\bf E} \left\{ X
      e^{-{c}X\, \beta}
  \right\}.
\end{equation}
\end{theo}

\bigskip
Observe the following duality of this result. For any $c<1$ we know
that the subgraph induced in our model by the long-range edges may
have at most $O(\log N)$ vertices in a connected component. According
to Theorem \ref{T}, for any $c<1$ there is 
\[q^{cr}(c)= \frac{1-c}{1+c}\]
such that for all $q^{cr}(c)<q<1$ our model will have a giant
component with a size of order $N$, while any  $q<q^{cr}(c)$ is
insufficient to produce  a giant component in $G_N(q,c)$ . Hence,
Theorem \ref{T} may also tell us something  about the "distances"
between the components of a random graph when it is considered on the
vertices of ${\bf Z}$.

\begin{rem}
In the proof of (\ref{E1}) we will show how to obtain $\alpha (q,c)$, and will discuss how optimal this value is. 
Statement (\ref{E2}) is rather trivial (and far from being optimal): 
it follows from a simple observation (see the details below) that
\[ {\bf P} \left\{ C_1  \Big( G_N(q,c)  \Big) <\, \frac{1}{|\log q|} \log N 
\right\}  \leq {\bf P} \left\{ C_1  \Big( G_N(q,0)  \Big) <\, \frac{1}{|\log q|} \log N 
\right\}  \rightarrow  0.\]
 
\end{rem}

\begin{rem}
For any fixed $c$ function $\beta (q,c)$ is continuous at $q=0$:
 if  $q=0$, i.e., when our graph is merely a classical $G_{n,c/n}$
  random graph, then $X \equiv 1$ and  (\ref{be}) becomes a
  well-known relation. Equation (\ref{be}) can be written in an exact form:
\[\beta =1-\frac{e^{c\beta}}{\left(e^{c\beta}-q\right)^2} (1-q)^2\, .\]
\end{rem}

  \bigskip
   It is easy to check that if $c \leq c^{cr}$ then the equation \eqref{be} does not have a
  strictly  positive
  solution, while $\beta=0$ is always a solution to \eqref{be}. 
Therefore one can  derive
  \begin{equation}\label{be'}
  \beta '_{c} \mid _{\ c \downarrow c^{cr} }=2 \frac{\big({\bf E}
    (X^2)\big)^3}{\big({\bf E} X\big)^2{\bf E} (X^3)}=2 \frac{(1-q^2)^2}{q^2+4q+1}.
\end{equation}
This shows that the emergence of the giant component at critical
parameter $c=c^{cr} $ becomes slower as $q$ increases, but the phase
transition remains  of the second order (exponent 1) for any $q<1$.

We
conjecture that similar results hold in the higher dimensions if
$q<Q^{cr}(d)$, where $Q^{cr}(d)$ is the percolation threshold in the
dimension $d$.
More exactly,
Theorem \ref{T} (as well as the first
equality in (\ref{be'})) should hold with $X$ replaced
by another random variable, which is stochastically not larger than
the size of the open  cluster at the origin in the edge percolation
model with a probability of edge $q$.
It is known from the percolation
theory (see, e.g., \cite{G}) that the tail of the
distribution of the size of an open cluster in the subcritical phase
decays exponentially. This should make possible to extend our arguments
(where we use essentially the distribution of $X$) to the general
case.

Our result in the supercritical case, namely  equation (\ref{be}) looks 
somewhat similar to the equation obtained in \cite{CL} for the "volume" (the sum of
degrees of the involved vertices)
of the giant component in the graph with a given sequence of the
expected degrees.
Note, however, that the model in
\cite{CL} (as well as the
derivations of the results)
differs essentially  from the one studied here.
In particular, in our model
 the critical mean degree 
when $c=c^{cr}$ and $N\rightarrow \infty$  is given by
\begin{equation}\label{deg}
2q+c^{cr} = 2q + \frac{1-q}{1+q} =
1 +  \frac{2q^2}{1+q} 
  \end{equation}
which  is strictly greater than 1 for all positive $q$.
This is in a contrast with the model studied in \cite{CL},
where the
critical expected average degree is still 1 as  in the classical random graph.

Although our model (when considered on the ring or torus in higher
dimensions)
is a perfectly {\it homogeneous} random graph, in the
sense that the degree distribution is the same for any  vertex, we
study it via {\it   inhomogeneous } random graphs, making use of the recently
developed theory from \cite{BJR}. The idea is the following. First, we consider  the
subgraph induced by the short-range edges, i.e., the edges which connect two neighbouring nodes with probability $q$. It is composed of the
consecutive connected paths (which may consist  just of one single vertex) on
$V_N=\{1,\ldots, N\}$. Call a
{\it macro-vertex} each of the component of this subgraph. We say that
a macro-vertex is of type $k$, if $k$ is the number of
vertices in it. 
Conditionally on the set of  macro-vertices, we
consider a graph on these macro-vertices induced by the long-range
connections. Two macro-vertices are said to be connected if there is
at least one (long-range type) edge  between two vertices
belonging to different macro-vertices.  
Thus the probability of an edge between two macro-vertices $v_i$ and $v_j$
of types $x$
and $y$ correspondingly,  is 
\begin{equation}\label{p}
    {\widetilde p}_{xy}(N):= 1-\left(1-\frac{c}{N} \right)^{x y}.
\end{equation}
Below we argue that this model fits the conditions of a general
inhomogeneous graph model defined in \cite {BJR},
find the critical parameters
and characteristics for the
graph on macro-vertices,  and then we turn back to the original model. 
We use essentially the results from \cite {BJR} to derive
(\ref{be}), while in the subcritical case our approach somewhat differs from the one in \cite {BJR}; we discuss this in the end of Section 2.4. We shall also note that our graph on macro-vertices is similar to the model studied in \cite{MR}, and our results on the
critical value agree
with those in \cite{MR}. 

Finally we comment that
our result should help to study
more general  model for the propagation of the neuronal activity
introduced in \cite{TV}. Here we show that
 a giant component in the graph can emerge from two sources, none of
 which can be neglected, but each of
 which may be in the subcritical phase, i.e., even  when both 
 $q<1$ and
 $c<1$. In particular, for any $0<c<1$ we can find $q<1$ which allows
 with a positive probability the
 propagation of  impulses  through the large part of the network due
 to the local activity.

\section{Proof}
\setcounter{equation}{0}

\subsection{Random graph on macro-vertices.}

Denote $X$ a random number of the vertices connected through
short-range edges to the vertex $1$ on $V_{\infty}^1 =\{1, 2,\ldots
\}$. Clearly, 
$X$ has the First success distribution defined in (\ref{X}).
  Let $X_1, X_2, \ldots $ , be independent copies of $X$, and define for any $N>1$
  \[T(N):=\min\{n\geq 0: \sum _{i=1}^n X_i \leq N, \ \ \sum _{i=1}^{n+1} X_i
  > N \},\]
where we assume that a sum over an empty set equals zero.
  
  Consider  now the
subgraph on $V_N=\{1, \ldots, N\}$ induced by the short-range
edges. This means that any two vertices $i$ and $i+1$ from $V_N$ are
connected with probability $q$  independent of the rest. By the
construction
this subgraph, call it $G_N^{(s)}(q)$,
is composed of
a random number of
connected paths of random  sizes. We call here the size of a path the
number of its vertices. Clearly, there is  a probability space
$(\Omega, {\cal F}, {\bf P})$ where the number of paths in $G_N^{(s)}(q)$
equals $T(N)$ if $\sum _{i=1}^{T(N)} X_i = N$, or
$T(N)+1$ if $\sum _{i=1}^{T(N)} X_i < N$. Correspondingly, the sizes
of the paths
follow the distribution of
\begin{equation}\label{t} 
{\bf X} = (X_1, X_2, \ldots , X_{T(N)}, N-\sum _{i=1}^ {T(N)}X_i)
\end{equation}
(where the last entry may take zero value). 

On the other hand, 
 the number of the connected components of $G_N^{(s)}(q)$ exceeds exactly by one the number of 
"missed" short edges on $V_N$. This means that on the same probability space $(\Omega, {\cal F}, {\bf P})$ there is a random variable $Y_N$ distributed as 
$Bin(N-1,1-q)$, such that either $T(N)=Y_N+1$ or $T(N)+1=Y_N+1$, and in any case
\begin{equation}\label{Y} 
0 \leq T(N)-Y_N \leq 1.
\end{equation}
This together with the Strong Law of Large Numbers implies
\begin{prop}\label{P1}
$$\frac{T(N)}{N} \ \stackrel{a.s.}{\rightarrow} \  1-q = \frac{1}{{\bf E}X} $$ 
as $N \rightarrow \infty$.$\hfill{\Box}$
\end{prop}

Also the relation (\ref{Y}) allows us to use the large deviation inequality 
from \cite{JLR} (formula (2.9), p.27 in \cite{JLR}) for the binomial random variables in order to obtain the following rate of convergence
\begin{equation}\label{T(N)} 
 {\bf P}\left\{ \left| \frac{T(N)}{N}- \frac{1}{{\bf E}X}
\right| > \delta \right\}
  \leq 2 \exp \left( -\frac{\delta^2}{12(1-q)} \, N\right)
\end{equation}
for all $\delta >0$ and $N>2/\delta$.

Define for any $k \geq 1$ an indicator function 
\begin{displaymath}
    I_{k}(x) = \left\{ \begin{array}{ll}
            1, & \mbox{ if } x=k, \\
            0, & \mbox{ otherwise. }
        \end{array} \right.
\end{displaymath}
As an immediate corollary of Proposition \ref{P1} and the Law of Large Numbers
we also get the following result.

\begin{prop}\label{P2} For any fixed $k\geq 1$
\begin{equation}\label{mu} 
  \frac{1}{T(N)} \sum _{i=1}^ {T(N)}  I_{k}(X_i) \
  \stackrel{P}{\rightarrow} \  {\bf P}
  \{X=k\}=(1-q)q^{k-1}=: \mu  (k)
\end{equation}
  as $N \rightarrow \infty$.$\hfill{\Box}$
\end{prop}

Given a  vector of paths ${\bf X}$ defined in (\ref{t}),   we
introduce another graph $ {\widetilde G}_{N}({\bf X}, q,c)$ as follows.
The set of vertices of $ {\widetilde G}_{N}({\bf X},q,c)$ we denote
$\{v_1, \ldots , v_{T(N)}\}$.  Each vertex $v_i$
 is
said to be of 
type $X_i$,
which means that 
 $v_i$ corresponds to the set of $X_i$ connected vertices on
 $V_N$. We shall also call any vertex $v_i$ of $ {\widetilde G}_{N}({\bf X}, q,c)$
a {\it macro-vertex}, and  write
\begin{equation}\label{v}
  v_i=\left\{
  \begin{array}{ll}
    \{1, \ldots, X_1\}, & \mbox{ if } i=1;\\ \\
    \{\sum_{j=1}^{i-1}X_j+1,
  \ldots, \sum_{j=1}^{i-1}X_j + X_i\}, & \mbox{ if }   i>1.
  \end{array}
  \right.
\end{equation}
With this notation the type of a vertex $v_i$ is simply the cardinality of set
$v_i$.
The space of 
 the types of macro-vertices
 is  $S=\{1,2, \ldots  \}$. According to (\ref{mu}) the
 distribution of type of a (macro-)vertex in graph $ {\widetilde
   G}_{N}({\bf X},q,c)$
 converges to measure $\mu$ on $S$.
  The edges between the vertices of $ {\widetilde G}_{N}({\bf X},q,c)$ are
 presented independently
 with probabilities induced by the original graph
 $G_{N}(q,c)$. More precisely,
the probability of an edge between any two vertices $v_i$ and $v_j$
of types $x$
and $y$ correspondingly,  is 
 ${\widetilde p}_{xy}(N)$
 introduced in (\ref{p}).
Clearly, this construction provides a one-to-one correspondence
between the connected components in the graphs $ {\widetilde
  G}_{N}({\bf X},q,c)$ and $
  G_{N}(q,c)$: the number
  of the connected components is the same
  for
  both
  graphs, as well as the number of the involved vertices from $V_N$
in two corresponding components. In fact, considering conditionally on
${\bf X}$ graph 
 $ {\widetilde
  G}_{N}({\bf X},q,c)$ we neglect only those long-range edges from 
$G_{N}(q,c)$, which connect vertices within each $v_i$, i.e., the
vertices which  are already
connected
 through the short-range edges.

  Consider now
\begin{equation}\label{p3}
    {\widetilde p}_{xy}( N) = 1-\left(1-\frac{c}{N} \right)^{xy}
    =: \frac{\kappa ' _N (x,y)}{N}.
\end{equation}
Observe that if $x(N)\rightarrow x$ and $y(N)\rightarrow y$ then
\begin{equation}\label{ka}
    \kappa ' _N (x(N), y(N))\rightarrow  c  xy
\end{equation}
 for all $x, y \in S$.
In order to place our
  model into the framework of the inhomogeneous random graphs
  from \cite{BJR} let us introduce another (random) kernel
  \[
    \kappa _{T(N)} (x,y)=\frac{T(N)}{N}  \kappa_N '(x, y),
\]
so that we can rewrite the probability ${\widetilde p}_{xy}(N)$
  in a graph
  $ {\widetilde
  G}_{N}({\bf X},q,c)$ taking into account the size of the graph:
\begin{equation}\label{pij}
    {\widetilde p}_{xy}( N) = \frac{\kappa _{T(N)} (x, y)}{T(N)}.
\end{equation}
(We use notations from \cite{BJR} whenever it is appropriate.)
According to Proposition \ref{P1} and (\ref{ka}),
if $x(N)\rightarrow x$ and $y(N)\rightarrow y$ then
\begin{equation}\label{kap}
    \kappa _{T(N)} (x(N), y(N))\rightarrow  \kappa (x,
    y):=\frac{c}{{\bf E}X}  xy \ \ \
    a.s.
\end{equation}
as $N \rightarrow \infty$ for all $x, y \in S$ .

Hence,  
in view of Proposition \ref{P2}
we conclude that
conditionally on $T(N)=t(N)$, where $t(N)/N \rightarrow 1/{\bf E}X$,
 our model 
falls into  the so-called "rank 1 case" of the general inhomogeneous
random graph model $G^{\cal V}(t(N),\kappa _{t(N)})$ with a vertex space
${\cal V}=(S,\mu, (X_1, \ldots, X_{t(N)})_{N\geq 1})$ from \cite{BJR} 
(Chapter 16.4). Furthermore, it is not difficult  to verify with a
help of the Propositions \ref{P1} and \ref{P2} that
\begin{equation}\label{c1}
\kappa \in L^1(S \times S, \mu \times \mu),
\end{equation}
since
\begin{equation*}
\sum_{y=1}^{\infty}\sum_{x=1}^{\infty} (1-q)x q^{x-1}(1-q)y q^{y-1} =
\Big( \frac{1}{1-q}\Big)^2,
\end{equation*}
and for any $t(N)$ such that $t(N)/N \rightarrow 1/{\bf E}X$
\begin{equation}\label{cond}
\frac{1}{t(N)} {\bf E}\{ e({\widetilde G}_{N}({\bf X},q,c))|T(N)=t(N)\}
\rightarrow
\frac{1}{2}\sum_{y=1}^{\infty} \sum_{x=1}^{\infty}\kappa (x,y) \mu(x)
\mu (y),
\end{equation}
where $e(G)$ denotes the number of edges in a graph $G$.
According to
Definition 2.7 from \cite{BJR},  under 
the conditions (\ref{cond}), (\ref{c1}) and (\ref{kap})
the sequence of  kernels $\kappa_{t(N)}$
(on the countable space $S \times S$)
is called {\it graphical} on
${\cal V}$
with limit $\kappa$.

\subsection{A branching process related to $ {\widetilde
  G}_{N}({\bf X},q,c)$.}
Here we closely follow the approach from \cite{BJR}.
We shall use a well-known technique of branching processes to reveal the
connected component in graph $ {\widetilde
  G}_{N}({\bf X},q,c)$. Recall first the usual algorithm of finding a
connected component. Conditionally on the set of macro-vertices, take any vertex $v_i$ to be the root. Find all
the vertices $\{ v^1_{i_1},v^1_{i_2},...,v^1_{i_n}\}$ connected to this vertex $v_i$ in the graph $ {\widetilde
  G}_{N}({\bf X},q,c)$, call them the first generation of $v_i$, and
then mark $v_i$ as "saturated". Then for each
non-saturated but already revealed vertex, we find all the vertices
connected to them but which have not been used previously. 
We continue this process until we end
up with a tree of saturated vertices.

Denote  $\tau_{N}(x)$ the set of the macro-vertices in
the tree constructed according to the  above algorithm with the
root at a vertex of type $x$.

It is plausible to think (and in our case it is correct, as will be
seen below) that
this algorithm with a high
probability as $N\rightarrow \infty$ reveals a tree of the offspring
of the
following multi-type 
 Galton-Watson process with type space $S=\{1,2, \ldots\}$: at any step,
 a particle of type $x \in S$ is replaced in the next
 generation by a set of particles where the number of  particles of
 type $y$ has a Poisson distribution $Po(\kappa
 (x,y) \mu(y))$.
 Let $\rho(\kappa;x)$ denote the probability that a particle of type $x$
 produces an infinite population.

\begin{prop}\label{P4}  The function $\rho(\kappa;x)$, $x \in S$, is the maximum
  solution to
\begin{equation}\label{rho}
\rho(\kappa;x) = 1- e^{-\sum_{y=1}^{\infty} \kappa
 (x,y) \mu(y)\rho(\kappa;y) }.
\end{equation}
  \end{prop}

\noindent  
{\bf Proof.}
We have
\begin{equation*}
\sum_{y=1}^{\infty}\kappa(x,y)\mu(y)= \frac{c}{{\bf E}X} \frac{x}{1-q} < \infty \text{ for
any $x$},
\end{equation*}
which together with (\ref{c1})
verifies that the conditions of
Theorem 6.1 from \cite{BJR} are satisfied, and the result (\ref{rho})
follows by
this theorem. $\hfill{\Box}$

\bigskip

Notice that it also follows by the same Theorem 6.1 from \cite{BJR}
that $\rho(\kappa;x)>0$ for all $x \in S$ if and only if
\begin{equation}\label{cr12}
  \frac{c}{{\bf E}X} \sum_{y=1}^{\infty}y^2\mu(y)=
  c \, \frac{{\bf E}X^2}{{\bf E}X}  =  c \, \frac{1+q}{1-q}>1;
\end{equation}
otherwise, $\rho(\kappa;x)= 0$ for all $x \in S$. Hence, the formula  (\ref{cr1})  for the
critical value follows from (\ref{cr12}).

As we showed above, conditionally on $T(N)$ (so that $T(N)/N \rightarrow 1/{\bf E}X$) the
sequence $\kappa_{T(N)}$ is graphical on $\cal
V$. Hence, the conditions of Theorem 3.1 from \cite{BJR} are satisfied
and we derive (first,
conditionally on $T(N)$, and therefore unconditionally)
that 
\[
\frac{C_1({\widetilde G}_{N}({\bf X}, q,c))}{T(N)} \stackrel{P}{\rightarrow} \rho(\kappa),
\]
where $\rho(\kappa)=\sum_{x =1}^{\infty}\rho(\kappa;x)\mu(x)$. This
together with Proposition \ref{P1} on the $a.s.$ convergence of $T(N)$
implies
\begin{equation}\label{MG}
\frac{C_1({\widetilde G}_{N}({\bf X}, q,c))}{N} \stackrel{P}{\rightarrow} (1-q)\rho(\kappa).
\end{equation}
Notice that here $C_1({\widetilde G}_{N}({\bf X}, q,c))$ is the number
of macro-vertices in ${\widetilde G}_{N}({\bf X}, q,c)$.

\subsection{On the distribution of types of vertices in $ {\widetilde
  G}_{N}({\bf X},q,c)$.}

 Given a  vector of paths ${\bf X}$ (see (\ref{t}))
we define a random sequence
\[{\bf \cal N} =\{{\cal N}_1, \ldots {\cal N}_N\} ,\]
where
\[
 {\cal N} _k= {\cal N} _k ({\bf X})= \sum _{i=1}^ {T(N)}  I_{k}(X_i). \]
In words,  ${\cal N} _k$ is  the number of (macro-)vertices of type $k$ in the set of
vertices of graph $ {\widetilde
  G}_{N}({\bf X},q,c)$. We shall prove here a useful result on the
distribution of ${\cal N}$ (which is stronger than Proposition \ref{P2}).

\begin{lem}\label{L1}
 For any fixed $\varepsilon >0$
\begin{equation}\label{x} 
{\bf P}
  \{  | {\cal N}_k/T(N) - \mu  (k)| > \varepsilon \, k \mu  (k) \ \
  \mbox{ for some } 1\leq k \leq N \}
  =o(1)
\end{equation}
as $N \rightarrow \infty$.
\end{lem}
\noindent
{\bf Proof.} Let us fix $\varepsilon >0$ arbitrarily. Observe that for any $K>1/\varepsilon$
\begin{equation}\label{A1}
{\bf P}
  \{  | {\cal N}_k/T(N) - \mu  (k)| > \varepsilon \, k \mu  (k) \ \
  \mbox{ for some } 1\leq k \leq N \}
\end{equation}
  \[ \leq {\bf P} \{ \max_{1\leq i \leq T(N)} X_i >
  K  \} 
\]
\[+ {\bf P}
  \{  | {\cal N}_k/T(N) - \mu  (k)| > \varepsilon \, k \mu  (k) \ \
  \mbox{ for some } 1\leq k \leq K \}.
\]
Next we shall choose an appropriate  $K=K(N)$ so that we will be able to
bound from above by $o(1)$ (as $N\rightarrow \infty$)
each of the summands on the right in (\ref{A1}).

Let us fix $\delta >0$ arbitrarily, and define an event
\begin{equation}\label{calA}
{\cal A}_{\delta , N} =\left \{ \left| \frac{T(N)}{N}- \frac{1}{{\bf E}X}\right|
\leq
\delta \right\}.
\end{equation}
Recall that according to (\ref{T(N)})
\begin{equation}\label{AA}
{\bf P} ({\cal A}_{\delta , N}) \geq 1- 2
\exp \left( -\frac{\delta^2}{12(1-q)} \, N\right) =1-o(1)
\end{equation}
as $N \rightarrow \infty$. Now we derive
\begin{equation}\label{A2}
{\bf P} \{  \max_{1\leq i \leq T(N)} X_i > K \}
\leq
  {\bf P} \{ 
    \max_{1\leq i \leq T(N)} X_i > K 
  \mid   {\cal A}_{\delta , N} \}{\bf P} ({\cal A}_{\delta , N}) +{\bf P}\{\overline{{\cal A}_{\delta, N}}\}
\end{equation}
\[ \leq
  \left(\frac{1}{{\bf E}X} +\delta \right)N \, {\bf P} \{ X_1 > K 
  \mid  {\cal A}_{\delta , N}  \}{\bf P} ({\cal A}_{\delta , N})
  +{\bf P}\{\overline{{\cal A}_{\delta, N}}\}
\]
\[ \leq
  \left(\frac{1}{{\bf E}X} +\delta \right)N \,
{\bf P} \{ X_1 > K \} 
 +{\bf P}\{\overline{{\cal A}_{\delta, N}}\}
\]
as $N \rightarrow \infty$.
Making use of the formula (\ref{X}) for the distribution of $X_1$
we obtain from (\ref{A2}) and (\ref{AA})
\begin{equation}\label{A3}
{\bf P} \{ \max_{1\leq i \leq T(N)} \geq K\} \leq
   C N q^K +2
\exp \left( -\frac{\delta^2}{12(1-q)} \, N\right) 
\end{equation}
as $N \rightarrow \infty$, where $C=C(\delta, q)$ is some finite
positive constant. Let now $\omega_1 (N)<N$
be any function tending to infinity with $N$, and set
\begin{equation}\label{A4}
K(N)= \frac{1}{|\log q|} \log N + \omega_1 (N).
\end{equation}
Clearly, bound (\ref{A3}) with $K$ replaced by $K(N)$ implies
\begin{equation}\label{A5}
{\bf P} \{ \max_{1\leq i \leq T(N)} X_i \geq K(N) \} =o(1)
\end{equation}
as $N \rightarrow \infty$.

Now we consider the last term in (\ref{A1}).
Let us define 
\begin{equation}\label{k0}
k_0:= \max \left\{ \left[ \frac{1}{\varepsilon}\right], \left[\frac{1}{|\log q|} \right] \right\} +2.
\end{equation}
Then we obtain
 making use of (\ref{AA})
\begin{equation}\label{A6}
{\bf P}
  \{  | {\cal N}_k/T(N) - \mu  (k)| > \varepsilon \, k \mu  (k) \ \
  \mbox{ for some } 1\leq k \leq K(N) \}
  \end{equation}
\[\leq
\sum_{k=1}^{k_0}
{\bf P}
  \{
  | \frac{1}{T(N)}{
    \sum _{i=1}^ {T(N)}  I_{k}(X_i)
  }- \mu  (k)|
  > \varepsilon \, k \mu  (k)
\} \]
\[
+
\sum_{k=k_0+1}^{K(N)}
{\bf P}
  \{
  | \frac{1}{T(N)}{
    \sum _{i=1}^ {T(N)}  I_{k}(X_i)
  }- \mu  (k)|
  > \varepsilon \, k \mu  (k)
\mid  {\cal A}_{\delta , N}
\}{\bf P} ({\cal A}_{\delta , N}) +o(1)
\]
\[=\sum_{k=k_0+1}^{K(N)}
{\bf P}
  \{
  | \frac{1}{T(N)}{
    \sum _{i=1}^ {T(N)}  I_{k}(X_i)
  }- \mu  (k)|
  > \varepsilon \, k \mu  (k)
\mid  {\cal A}_{\delta , N}
\}{\bf P} ({\cal A}_{\delta , N}) +o(1)\]
as $N\rightarrow \infty$,
where the last equality is due to Proposition \ref{P2}. Notice that 
for each $k>k_0$ we have $\varepsilon \, k >1$ and therefore
\begin{equation}\label{A7}
{\bf P}
\left\{
\left| \frac{1}{T(N)}{
  \sum _{i=1}^ {T(N)}  I_{k}(X_i)}- \mu  (k)\right|
  > \varepsilon \, k \mu  (k)
\mid  {\cal A}_{\delta , N}
\right\}
\end{equation}
\[={\bf P}
\left\{
 \frac{1}{T(N)}{
  \sum _{i=1}^ {T(N)}  I_{k}(X_i)}- \mu  (k)
  > \varepsilon \, k \mu  (k)
\mid  {\cal A}_{\delta , N}
\right\}
\]
\[ \leq {\bf P}
\left\{
\frac{1}{\left(\frac{1}{{\bf E}X} -\delta \right)N}
{\sum _{i=1}^{\left[\left(\frac{1}{{\bf E}X} +\delta \right)N\right]+1}
    I_{k}(X_i)}
  > \mu  (k) + \varepsilon \, k \mu  (k)
\mid  {\cal A}_{\delta , N}
\right\}
=:{\bf P}(k).
\]
 Set $t(N)=\left[
        \left(\frac{1}{{\bf E}X} +\delta \right)N
      \right]+1$. Then 
using the bound
\[
\frac{
\left(
\frac{1}{{\bf E}X} -\delta 
\right)N 
}{t(N)} > 1- \frac{5}{2} 
{\bf E}X \, \delta \]
for all $N>2/{\delta}$, we derive
\begin{equation}\label{A8}
{\bf P}(k) \leq
{\bf P}
\left\{
\frac{1}{t(N)
    }
\sum _{i=1}^{t(N)}
    I_{k}(X_i)
  > \mu  (k)(1+  \varepsilon \, k )(1- \frac{5 {\bf E}X}{2}\delta)
\mid  {\cal A}_{\delta , N}
\right\}
\end{equation}
 for all $N>2/{\delta}$. 
Now for all $k>k_0$ and 
$
0<\delta <\frac{1}{10{\bf E}X}$
we have $(1+  \varepsilon \, k )(1-\frac{5 {\bf E}X}{2} \delta)\geq 
1+  \frac{\varepsilon}{2} \, k$, and therefore
\begin{equation}\label{A9}
{\bf P}(k) \leq
{\bf P}
\left\{
\frac{1}{t(N)
    }
\sum _{i=1}^{t(N)}
    I_{k}(X_i)
  > \mu  (k)(1+  \frac{\varepsilon}{2} \, k )
\mid  {\cal A}_{\delta , N}
\right\}
\end{equation}
\[\leq
{\bf P}
\left\{
\frac{1}{t(N)
    }
\sum _{i=1}^{t(N)}
    I_{k}(X_i)
  > \mu  (k)(1+  \frac{\varepsilon}{2} \, k )
\right\} \, / \, {\bf P}({\cal A}_{\delta , N}). \]
Note that $\sum _{i=1}^{t(N)}
    I_{k}(X_i)$ follows the binomial distribution
$Bin(t(N),\mu (k))$. This allows us to use the large deviation inequality 
from \cite{JLR} (see (2.5), p.26 in \cite{JLR}) and derive
\begin{equation}\label{A10}
{\bf P}
\left\{
\frac{1}{t(N)
    }
\sum _{i=1}^{t(N)}
    I_{k}(X_i)
  > \mu  (k)(1+  \frac{\varepsilon}{2} \, k )
\right\}
\end{equation}
\[\leq \exp \bigg( - \frac{(\frac{\varepsilon}{2}  k \mu (k)
    t(N))^2}{\frac{1}{3} \varepsilon k \mu (k) t(N) + 2\mu (k)
    t(N)}\bigg)
 \leq \exp \big(- \frac{1}{10} \varepsilon k \mu (k) t(N)\big)
\]
for all $k>k_0$. Substituting this into (\ref{A9}) we obtain
\begin{equation}\label{A11}
{\bf P}(k) \leq  \exp \big(- \frac{1}{10} \varepsilon k \mu (k) t(N)\big)
\, / \, {\bf P}({\cal A}_{\delta , N})
\end{equation}
for all $k>k_0$. The last bound combined with (\ref{A7}) and (\ref{A6}) leads to
\[
{\bf P}
  \{  | {\cal N}_k/T(N) - \mu  (k)| > \varepsilon \, k \mu  (k) \ \
  \mbox{ for some } 1\leq k \leq K(N) \}
\]
\[\leq
\sum_{k=k_0+1}^{K(N)}\exp \big(- \frac{1}{10} \varepsilon k \mu (k) t(N)\big)
\, 
+o(1),\]
as $N\rightarrow \infty$. Taking into account  that function $k \mu (k)$
is decreasing for $k>k_0$ we derive from the last bound:
\begin{equation}\label{A13}
{\bf P}
  \{  | {\cal N}_k/T(N) - \mu  (k)| > \varepsilon \, k \mu  (k) \ \
  \mbox{ for some } 1\leq k \leq K(N) \}
  \end{equation}
\[\leq 
\exp \big(- \frac{1}{10} \varepsilon K(N) \mu (K(N)) t(N) + \log K(N)\big)
+o(1),\]
as $N\rightarrow \infty$.

Setting now $\omega_1 (N) = \log \log \log N$ in (\ref{A4}), it is easy to check that for
$$K(N) = \frac{1}{|\log q|} \log N + \log \log \log N $$
the entire right-hand side of the inequality (\ref{A13}) is $o(1)$
 as $N\rightarrow \infty$. This together with the previous bound (\ref{A5}) and inequality (\ref{A1})
finishes the proof of lemma. \hfill$\Box$

\subsection{Proof of Theorem \ref{T} in the subcritical case $c  < c^{cr} (q)$  .}
Let us fix $0\leq q < 1$ and then $c  < c^{cr} (q)$  arbitrarily.
Given ${\bf X}$ let again $v_i$ denote the macro-vertices with types
$X_i$, $i=1,2, \ldots,$
respectively, and  let ${\widetilde
  L}$
denote a connected component in 
${\widetilde
  G}_{N}({\bf X},q,c)$. Firstly, for any $K>0$ and $0<\delta<1/{\bf E}X$
we derive with help of (\ref{AA})
\begin{equation}\label{S5}
{\bf P} \left\{  C_1\Big( G_N(q,c)  \Big) < K
\right\}
\leq {\bf P} \left\{  C_1\Big( G_N(q,0)  \Big) < K
\right\} =
{\bf P} \left\{ \max _{1\leq i \leq T(N)} X_i
<K \right\} 
\end{equation}
\[\leq {\bf P} \left\{ \max _{1\leq i \leq T(N)} X_i
<K \mid {\cal A}_{\delta, N}\right\} +o(1) \leq \left(1-{\bf P} \left\{ X
\geq K \right\} \right)^{N\left(\frac{1}{{\bf E}X}-\delta \right)} +o(1),\]
as $N\rightarrow \infty$, where $X$ has the
 $Fs(1-q)$-distribution.
Since 
\[{\bf P} \left\{ X
\geq K \right\} = q^{K-1},\]
we derive from (\ref{S5}) for any $a_1 < \frac{1}{|\log q|}$ and $K=a_1 \log N$
\[
{\bf P} \left\{  C_1\Big( G_N(q,c)  \Big) < a_1 \log N
\right\}
= o(1),\]
which proves statement (\ref{E2}).

Consider now
for any positive constant $a$ and a function $w=w(N)\geq \log N$
\begin{equation}\label{J1}
{\bf P} \left\{  C_1\Big( G_N(q,c)  \Big) > aw 
\right\}
=
{\bf P} \left\{ \max _{{\widetilde
  L}} \sum _{v_i \in {\widetilde
  L}}X_i
> a w \right\} .
\end{equation}
We know already from (\ref{MG})
that in the subcritical case
the
size (the number of macro-vertices) of any ${\widetilde
  L}$ is {\bf whp}  $o(N)$.
Note that  when the kernel $\kappa(x,y)$
is not bounded uniformly in both arguments, which is our case,
it is not granted that the
largest component in the subcritical case is at most of order $\log
N$ (see, e.g., discussion of 
Theorem 3.1 in \cite{BJR}). Therefore first we shall prove the following
intermediate result.

\begin{lem}\label{LS} If $c<c^{cr}(q)$ then 
  \begin{equation}\label{S11}
{\bf P} \left\{  C_1\Big( {\widetilde
  G}_N({\bf X}, q,c)  \Big) > N^{1/2}
\right\}
=o(1).
\end{equation}
\end{lem}

\noindent
{\bf Proof.}
Let us fix $\varepsilon >0$ and $\delta >0$ arbitrarily and introduce
the following event
\begin{equation}\label{A19}
{\cal B}_N:=  {\cal A}_{\delta, N}
 \cap \left( \max _{1\leq i\leq T(N)} X_i \leq \frac{2}{|\log q|}\log N \right) \cap \left( \cap _{k=1}^N \left\{\left|  \frac{{\cal N}_k}{T(N)} -\mu(k)\right|\leq \varepsilon k \mu(k)
\right\} \right).
\end{equation}
According to (\ref{AA}),  (\ref{A5}) and
(\ref{x}) we have
\begin{equation}\label{cb1}
{\bf P} \left\{ {\cal B}_N \right\} =1 -o(1)
\end{equation}
as $N \rightarrow \infty$.

Recall that  $\tau_{N}(x)$ denotes the set of the macro-vertices in
the tree constructed according to the algorithm of revealing of
connected component described above.
Let $|\tau_{N}(x)|$ denotes the
number of macro-vertices in
$\tau_{N}(x)$.
Then we easily derive
\begin{equation}\label{SA18}
 {\bf P} \left\{  C_1\Big( {\widetilde
  G}_N({\bf X}, q,c)  \Big) > N^{1/2}
\right\}
\leq
{\bf P} \left\{ \max _{1\leq i\leq T(N)} |\tau_{N}(X_i)| >  N^{1/2}\mid {\cal B}_N 
\right\} +o(1)
\end{equation}
\[
\leq
N \sum_{k=1}^{N} 
(1+ \varepsilon k) \mu(k) \Big( \delta + 1/{\bf E}X \Big)
{\bf P} \left\{  
|\tau_{N}(k)| >  N^{1/2}\mid {\cal B}_N 
\right\} +o(1)
\]
as $N \rightarrow \infty$.
We shall use the multi-type branching process introduced above
(Section 2.2)
to approximate the distribution of $|\tau_{N}(k)|$.
Let further ${\cal X}^{c,q}(k)$ denote the number of the particles
(including the initial one) in
the branching process
starting with a single particle
of  type $k$. 
Observe that at each step of the exploration algorithm,
the number of new neighbours of $x$ of type $y$ has a
binomial distribution $ Bin(N_y',{\widetilde  p}_{xy}(N))$ where $N_y'$
is the number of remaining vertices of type $y$, so that 
$ N_y'
\leq {\cal N}_y $. 

We shall explore the following obvious 
relation between the Poisson and the binomial distributions. 
Let $Y_{n,p} \in Bin(n,p)$ and $Z_{a} \in Po(a)$, where $0<p<1/4$ and $a>0$.
Then for all $k\geq 0$
\begin{equation}\label{A21}
{\bf P}
  \{ Y_{n,p}=k\}
\leq (1+Cp^2 )^n \, 
{\bf P}\{ Z_{n\frac{p}{1-p}}=k
  \},
\end{equation}
where $C$ is some positive constant (independent of $n$ and $p$).
Notice that for all $x,y \leq \frac{2}{|\log q|}\log N$ 
\begin{equation}\label{pB}
    {\widetilde p}_{xy}(N)= 1-\left(1-\frac{c}{N} \right)^{x y}
    =\frac{c}{N} \ {x y} \ (1+o(1)),
\end{equation}
and clearly,  ${\widetilde p}_{xy}(N) \leq 1/4$ for all large $N$.
 Therefore 
for any fixed 
positive $\varepsilon _1$
we can choose small $\varepsilon$ and $\delta$ in (\ref{A19}) so that
conditionally on ${\cal B}_N$ we have
\begin{equation}\label{pB1}
    N_y' \frac{{\widetilde p}_{xy}(N)}{1-{\widetilde p}_{xy}(N)}
 \leq  (1+ y \varepsilon _1)\mu (y)\kappa (x,y)
\end{equation}
for all large $N$. Let us write further
\[ \mu(y) =\mu _{q} (y), 
\ \ \   \  \mu _{q}= \sum_{y\geq 1} y \mu _{q} (y)\ (={\bf E}
X), \   \   \   \         \   \kappa (x,y)=   \kappa_{c,q} (x,y) \]
emphasizing dependence on $q$ and $c$.
Then for any 
$\varepsilon _2>0$ and any 
$q'>q$ such that $\frac{1-q}{q} \frac{q'}{1-q'}<1+\varepsilon_2$
 we can choose  $\varepsilon _1 < \log(q'/q)$, and derive from  (\ref{pB1})
\begin{equation}\label{pB2}
    N_y' \frac{{\widetilde p}_{xy}(N)}{1-{\widetilde p}_{xy}(N)}
 \leq  (1+\varepsilon _2)\mu_{q'}  (y)\kappa_{c,q}  (x,y)=
\mu_{q'}  (y) \frac{(1+\varepsilon _2) \, c}{\mu _{q} }xy .
\end{equation}
Setting now
$c':=(1+\varepsilon _2)\frac{\mu_{q'} }{\mu_{q} } \, c$
we rewrite (\ref{pB2}) as follows
\begin{equation}\label{pB3}
    N_y' \frac{{\widetilde p}_{xy}(N)}{1-{\widetilde p}_{xy}(N)}
 \leq  \mu_{q'}  (y)\kappa _{c', q'} (x,y).
\end{equation}
Recall that above  we fixed $q$ and $c<c^{cr}(q)$, where $c^{cr}(q)$ 
is decreasing and continuous in $q$. Hence, we can choose
 $q'>q$ and $c':=(1+\varepsilon _2)\frac{\mu_{q'} }{\mu_{q} }\, c $
so that 
\begin{equation}\label{A22}
    c<c'<c^{cr}(q')<c^{cr}(q),
\end{equation}
and moreover $c'$ and $q'$ can be chosen arbitrarily close to $c$ and $q$, respectively.

Now conditionally on ${\cal B}_N$
we can replace according to (\ref{A21})  at each (of at most $N$) step
of the exploration algorithm the
$ Bin(N_y',{\widetilde  p}_{xy}(N))$ variable 
with $Po(N_y' \frac{{\widetilde p}_{xy}(N)}{1-{\widetilde p}_{xy}(N)})$,
and further replace the last variables with the stochastically larger
ones $Po(\mu_{q'}  (y)\kappa_{c', q'}  (x,y))$ (recall (\ref{pB3})).
As a result we get the 
following bound using branching process:
\begin{equation}\label{S12}
{\bf P}
\left\{ |\tau_N(k) | > N^{1/2}\mid {\cal B}_N 
\right\}
\end{equation}
\[
\leq
\left(1+C \left(\max_{x,y \leq 2 \log N/|\log q|} {\widetilde p}_{xy}(N)
\right)^2\right)^{N^2}
 \, {\bf P} \left\{ {\cal X}^{c',q'} (k)> N^{1/2}
\right\} .
\]
 This together with (\ref{pB}) implies
\begin{equation}\label{SA23}
{\bf P} \left\{ |\tau_N(k) | > N^{1/2} \mid {\cal B}_N 
\right\} 
\leq
e^{b(\log N)^4} \, 
{\bf P} \left\{
  {\cal X}^{c',q'} (k)> N^{1/2}
\right\} ,
\end{equation}
 where
$b $ is some positive constant. Substituting the last bound into (\ref{SA18})
we derive
\begin{equation}\label{SA24}
{\bf P} \left\{   C_1\Big( {\widetilde
  G}_N({\bf X}, q,c)  \Big) > N^{1/2}\right\} 
\end{equation}
\[
\leq
b_2N e^{b(\log N)^4}
\sum_{k=1}^{N} 
 k \mu _{q'}(k)
{\bf P} \left\{ {\cal X}^{c',q'} (k)> N^{1/2}
\right\}
 +o(1)
\]
as $N \rightarrow \infty$, where
$b_2$ is some positive constant. By the Markov's inequality
\begin{equation}\label{SMark}
{\bf P} \big\{{\cal X}^{c',q'} (k) > N^{1/2}\big\} \leq z^{-N^{1/2}}
{\bf
    E} z^{{\cal X}^{c',q'} (k)}
\end{equation}
for all $z\geq 1$.
Denote $h_z(k)={\bf E} z^{{\cal X}^{c',q'} (k)}$; then with a help of 
  (\ref{SMark}) we get from (\ref{SA24})
\begin{equation}\label{SA25}
{\bf P} \left\{  C_1\Big( {\widetilde
  G}_N({\bf X}, q,c)  \Big) > N^{1/2}  \right\} 
\leq b_1N e^{b(\log N)^2} z^{-N^{1/2}}
\sum_{k=1}^{N} 
 k \mu _{q'}(k)
h_z(k) + o(1).
\end{equation}
Now we will show that there exists $z > 1$ such that the series 
$$B_z (c', q')= \sum_{k=1}^{\infty} k \mu_{q'}(k) h_z(k)$$ 
 converge. This together with (\ref{SA25}) will clearly imply the
 statement of the lemma.

Note that  
function $h_z(k)$ (as a generating function for a  branching process) satisfies the following equation
\[
\begin{array}{ll}
h_z(k) & =z \exp {\left\{
\sum_{x=1}^{\infty} \kappa_{c', q'}  (k,x) \mu_{q'}(x) (h_z(x)-1 )
\right\}} \\
\\
& =z\exp \left\{
\frac{c'}{\mu_{q'}} k \left(
\sum_{x=1}^{\infty} x \mu_{q'}(x)h_z(x)-\mu_{q'}\right)\right\} \\
\\
& =z \exp \left\{ \frac{c'}{ \mu_{q'}} k  
( B_z(c', q')-\mu_{q'})\right\}.
\end{array}
\]
Multiplying both sides by $k \mu_{q'} (k)$ and summing up over $k$ 
we find
\begin{equation*}
B_z(c', q') = \sum_{k=1}^{\infty} k \mu_{q'} (k)z \exp \left\{\frac{c'}{\mu_{q'}} 
 k ( B_z(c', q')-\mu_{q'})\right\} .
\end{equation*}
Let us write  for simplicity $B_z=B_z(c, q)$.
Hence, as long as $B_z$ is finite, it  should satisfy equation
\begin{equation}\label{Az}
B_z = z{\bf E} X \, e^{\frac{c}{ {\bf E} X } X ( B_z-{\bf E} X )},
\end{equation}
which implies in turn that  $B_z$ is finite for some $z>1$ if and only if (\ref{Az}) has at least
one solution (for the same value of $z$). 
Notice that
\begin{equation}\label{Sinit}
B_z \geq B_1 = {\bf E} X 
\end{equation}
for  $z \geq 1$. Let us fix $z> 1$ and consider equation
\begin{equation}\label{SAz1}
y/z  = {\bf E} X \, e^{\frac{c}{ {\bf E} X } X ( y-{\bf E} X )}=:F(y)
\end{equation}
for $ y \geq {\bf E} X $.
Using the properties of the distribution of $X$ it is easy to derive 
that
function $$F(y)
=\frac{1}{{\bf E}X} \frac{ e^{c(\frac{y}{{\bf E}X}-1)}}{\big( 1- q
  e^{c(\frac{y}{{\bf E}X}-1)}\big)^2}$$
is 
  increasing
  and has positive second derivative if ${\bf E} X \leq y\leq y_0
  $, where $y_0$ is the root of $1= q
  \exp \left\{c(\frac{y_0}{{\bf E}X}-1)\right\}$. Compute now
\begin{equation}\label{SAz2}
\frac{\partial }{\partial  y}F(y)|_{y={\bf
    E}X } = 
\frac{c}{ {\bf E} X } {\bf E} X^2 =\frac{c}{c^{cr}}.
\end{equation}
Hence,  if $c<c^{cr}$ then there exists $z>1$ such that 
 there is a finite solution $y$ to (\ref{SAz1}). 
 Taking into account
condition (\ref{A22}),   we find that  $
B_z(c', q')$ is also finite for some $z>1$, which 
 finishes the proof of the lemma.
\hfill$\Box$

\bigskip

Now we are ready to complete the proof of (\ref{E1}), following almost the same arguments as in the proof of the previous lemma.
Let $S_{N}(x)= \sum_{v_i \in \tau_{N}(x)}X_i$ denote the
number of vertices from $V_N$ which compose the macro-vertices of
$\tau_{N}(x)$.
Denote
\[ {\cal B}'_N := {\cal B}_N \cap
\left( 
\max _{1\leq i \leq T(N)} |\tau_{N}(X_i)|<N^{1/2}
\right).
\]
According to (\ref{cb1}) and Lemma \ref{LS}
we have
\[{\bf P} \left\{ {\cal B}'_N \right\}=1- o(1).\]
This allows us to  derive from (\ref{J1}) 
\begin{equation}\label{A18}
{\bf P} \left\{  C_1 \Big({
  G}_{N}(q,c) \Big) > a w \right\} \leq
{\bf P} \left\{ \max _{1\leq i\leq T(N)} S_N(X_i) > a w \mid
  {\cal B}'_N 
\right\} +o(1)
\end{equation}
\[
\leq
N \sum_{k=1}^{N} 
(1+ \varepsilon k) \mu(k) \Big( \delta + 1/{\bf E}X \Big)
{\bf P} \left\{  
S_N(k) > a w \mid {\cal B}'_N 
\right\} +o(1).
\]

Let now $S^{c,q}(y)$  denote the 
sum of types (including the one of the initial particle) in  the total progeny
of the introduced above branching process 
starting with initial particle of type $y$.
Repeating the same argument which led to (\ref{S12}),
we get the
following bound using the introduced branching process:
\[
{\bf P} \left\{ S_N(k) > a w \mid {\cal B}'_N 
\right\} 
\leq
\left(1+C \left(\max_{x,y \leq 2 \log N/|\log q|} {\widetilde p}_{xy}(N)
\right)^2\right)^{b_1 N \sqrt{N}}
 \, {\bf P} \left\{ S^{c',q'} (k)> a w 
\right\} 
\]
as $N \rightarrow \infty$, where we take into account that we can
perform at most $\sqrt{N}$ steps of exploration (the maximal possible
number of  macro-vertices in any ${\widetilde L}$).
This  together with (\ref{pB}) implies
\begin{equation}\label{A23}
{\bf P} \left\{ \tau_N(k) > a w \mid {\cal B}'_N 
\right\} 
\leq
(1+o(1))
{\bf P} \left\{ S^{c',q'} (k)> a w 
\right\} 
\end{equation}
as $N\rightarrow \infty$. Substituting the last bound into (\ref{A18})
we derive
\begin{equation}\label{A24}
{\bf P} \left\{  C_1 \Big({
  G}_{N}(q,c) \Big) > a w \right\} 
\leq
N b
\sum_{k=1}^{N} 
 k \mu _{q'}(k)
{\bf P} \left\{ S^{c',q'} (k)> a w
\right\}
 +o(1)
\end{equation}
as $N \rightarrow \infty$, where $b$ is some positive constant. 
Denote $g_z(k)={\bf E} z^{S^{c',q'} (k)}$; then 
similar to (\ref{SA25})
 we derive from (\ref{A24})
\begin{equation}\label{A25}
{\bf P} \left\{  C_1 \Big({
  G}_{N}(q,c) \Big) > a w(N) \right\} 
\leq
b
N \sum_{k=1}^{N} 
 k \mu _{q'}(k)
g_z(k) z^{-a w(N)} + o(1).
\end{equation}
We shall search for all  $z\geq 1$ for which the series 
$$A_z (c', q')= \sum_{k=1}^{\infty} k \mu_{q'}(k) g_z(k)$$ 
 converge.
Function $g_z(k)$ (as a generating function for a certain branching process) satisfies the following equation
\[
\begin{array}{ll}
g_z(k) & =z^k \exp {\left\{
\sum_{x=1}^{\infty} \kappa_{c', q'}  (k,x) \mu_{q'}(x) (g_z(x)-1 )
\right\}} \\
\\
& =z^k \exp \left\{
\frac{c'}{\mu_{q'}} k \left(
\sum_{x=1}^{\infty} x \mu_{q'}(x)g_z(x)-\mu_{q'}\right)\right\} \\
\\
& =z^k \exp \left\{ \frac{c'}{ \mu_{q'}} k  
( A_z(c', q')-\mu_{q'})\right\}.
\end{array}
\]
Multiplying both sides by $k \mu_{q'} (k)$ and summing up over $k$ 
we find
\begin{equation}\label{SAz}
A_z(c', q') = \sum_{k=1}^{\infty} k \mu_{q'} (k)z^k \exp \left\{\frac{c'}{\mu_{q'}} 
 k ( A_z(c', q')-\mu_{q'})\right\} .
\end{equation}
It follows from here (and the fact that $A_z(c', q')\geq
\mu_{q'}$ for all $z\geq 1$) that  if there exists $z>1$ for which the
series $A_z(c', q')$ converge, it should satisfy
\begin{equation}\label{con1}
z<\frac{1}{q'}.
\end{equation}
According to  (\ref{SAz}), as long as 
 $A_z=A_z(c, q)$ is finite it  satisfies  the equation
\[
A_z = {\bf E} X z^X \, e^{\frac{c}{ {\bf E} X } X ( A_z-{\bf E} X )},
\]
which implies that $A_z$ is finite for some $z>1$
if and only if the last equation
has  at least one solution
\begin{equation}\label{init}
A_z \geq A_1 = {\bf E} X.
\end{equation}
Let us fix $z> 1$ and consider equation
\begin{equation}\label{Az1}
y = {\bf E} X z^X \, e^{\frac{c}{ {\bf E} X } X ( y-{\bf E} X )}=:f(y,z).
\end{equation}
Using the properties of the distribution of $X$ it is easy to derive 
that
$$f(y,z)=\frac{1}{{\bf E}X} \frac{z e^{c(\frac{y}{{\bf E}X}-1)}}{\big( 1- qz
  e^{c(\frac{y}{{\bf E}X}-1)}\big)^2}. $$
We shall consider
$f(y,z)$ for  $y\geq {\bf E}X$ and  $1\leq z<
  \frac{1}{q}e^{-c(\frac{y}{{\bf E}X}-1)}$.
  It is easy to check that in this area 
function $f(y,z)$ is
  increasing, it has all the derivatives of the
second order, and $\frac{\partial ^2}{\partial y^2} f(y,z)>0$. Compute now
\begin{equation}\label{Az2}
\frac{\partial }{\partial  y}f(y,z)|_{y=1, z=1} = 
\frac{c}{ {\bf E} X } {\bf E} X^2 =\frac{c}{c^{cr}}.
\end{equation}
Hence, if 
$c>c^{cr}$ there is no solution $y\geq {\bf E}X$ to (\ref{Az1}) for any $z>1$.
On the other hand, if $c<c^{cr}$ then there exists $1<z_0<1/q$ such that for all
$1\leq z <z_0$ there is a finite solution $y\geq {\bf E}X$ to (\ref{Az1}). One could find $z_0$  for example, as the (unique!) value for which function $y$ is tangent to $f(y,z_0)$ if $y\geq {\bf E}X$.

Now 
taking into account that $c'>c$ and $q'>q$ can be chosen arbitrarily close to $c$ and $q$, respectively, we derive from (\ref{A25})
that for all $1<z<z_0$ 
\begin{equation}\label{S1}
{\bf P} \left\{  C_1 \Big({
  G}_{N}(q,c) \Big) > a w(N) \right\} 
\leq b(z)
N  z^{-a w(N)} + o(1)
\end{equation}
as $N\rightarrow \infty$, where $b(z)<\infty$. This implies that for any $a>1/\log z_0>1/|\log q|$
\begin{equation}\label{S2}
{\bf P} \left\{  C_1 \Big({
  G}_{N}(q,c) \Big) > a \log N \right\} 
= o(1)
\end{equation}
as $N\rightarrow \infty$, which proves (\ref{E1}).
\hfill$\Box$

\bigskip

To conclude this section we comment on the methods 
used here. It is shown in \cite{T3} that in the subcritical case of classical random graphs the same method of generating functions combined with the Markov inequality leads to a constant which is known to be the principal term for the asymptotics of the size of the largest component (scaled to $\log N$). 
This gives us hope that
a constant $a$ chosen here to satisfy 
$a>1/\log z_0$ is close to 
 the minimal constant for which
 statement (\ref{S2}) still holds. 

Similar methods were used in \cite {T2} for 
some class of inhomogeneous random graphs, and in \cite{BJR} for a general class of models.
Note, however, some difference with the approach in \cite{BJR}.
 It is assumed in \cite{BJR}, Section 12,
that the generating function
for the corresponding branching process with the initial state $k$ (e.g., our function $g_z(k)$, $k\geq 1$)
is bounded uniformly in $k$. 
As we prove here this condition is not always necessary: we need only convergence of the series $A_z$, while $g_z(k)$ is  unbounded in $k$ in our case. Furthermore, our approach allows one to construct constant $\alpha(q,c)$ as a function of the parameters of the model.

  \subsection{Proof of Theorem \ref{T} in the supercritical case.}

Let $\mathcal{C}_k$ denote the set of vertices in the
$k$-th largest component in graph ${G}_N(q,c)$, and conditionally on ${\bf X}$
let $\widetilde{\mathcal{C}}_k$ denote the set of macro-vertices in the
$k$-th largest component in graph $\widetilde{G}_N({\bf X}, q,c)$ (ordered in any way if there are ties).
Let also ${C}_k$ and $\widetilde{{C}}_k$  denote correspondingly, their sizes.
According to our construction  for any connected component $\widetilde{L}$ in  
 $\widetilde{G}_N({\bf X}, q,c)$ there is a unique component $L$ in
 ${G}_N(q,c)$ such that they are composed of the same vertices from
 $V_N$, i.e., in the notations (\ref{v})
\[L=\cup _{v \in \widetilde{L}}\cup _{k \in v}\{k\}=:V(\widetilde{L}).\]
Next we prove that with a high probability the largest components in both graphs consist of the same vertices.

\begin{lem} \label{P5}
For any $0\leq q<1$ if $c>c^{cr}(q)$ then 
\begin{equation}\label{C}
{\bf P}\{{\mathcal{C}}_1 =V(\widetilde{\mathcal{C}}_1)\} =1-o(1)
\end{equation}
as $N \rightarrow \infty$.
\end{lem}
\noindent
{\bf Proof.} In a view of the argument preceeding this lemma we have
\[{\bf P}\{{\mathcal{C}}_1 \neq V(\widetilde{\mathcal{C}}_1)\} 
={\bf P}\{{\mathcal{C}}_1 = V(\widetilde{\mathcal{C}}_k) \mbox{ for some } k\geq 2\} .\]
According to Theorem 12.6 from \cite{BJR}, 
conditions of which are satisfied here, 
in the supercritical case conditionally on $T(N)$ such that $T(N)/N
\rightarrow 1/{\bf E}X$,
we have {\bf whp} 
$\widetilde{C}_2=O(\log(T(N)))$, which by  Proposition
\ref{P1} implies $\widetilde{C}_2=O(\log N )$ {\bf whp} .
Also we know already from (\ref{MG})
that in the supercritical case
  $\widetilde{C}_1=O(N)$ {\bf whp}, and therefore
${C}_1=O(N)$ {\bf whp}. Hence, 
for some positive constants $a$ and $b$
\begin{equation}\label{C1}
{\bf P}\{{\mathcal{C}}_1 \neq V(\widetilde{\mathcal{C}}_1)\} 
={\bf P}\{{\mathcal{C}}_1 = V(\widetilde{\mathcal{C}}_k) \mbox{ for some } k\geq 2\} 
\end{equation}
\[ \leq
{\bf P} \left\{\left( 
\max _{k\geq 2}|V(\widetilde{\mathcal{C}}_k)| >b  N 
\right) \cap \left( \max _{k\geq 2}\widetilde{{C}}_k < a \log N \right)
\right\} 
+o(1).
\]
Define now for any $K=K(N)$ a set $$B_N:=\{\exists X_i \geq K \text{ for some } 1\leq i \leq T(N)\}.$$ According to (\ref{A3})
\[
{\bf P} (B_N) \leq
   C N q^K +o(1)
\]
as $N \rightarrow \infty$ for some constant $C$ independent of $K$ and $N$.
Setting from now on $K=\sqrt{N}$ we have
$
{\bf P} (B_N) = o(1)$
as  $N \rightarrow \infty$.
 Then we derive
\begin{equation}\label{C2}
{\bf P} \left\{\left( 
\max _{k\geq 2}|V(\widetilde{\mathcal{C}}_k)| >b N 
\right) \cap \left( \max _{k\geq 2}\widetilde{{C}}_k < a \log N \right)
\right\} 
\end{equation}
\[ \leq 
{\bf P} \left\{\left( 
\max _{k\geq 2}|V(\widetilde{\mathcal{C}}_k)| >b N 
\right) \cap \left( \max _{k\geq 2}\widetilde{{C}}_k < a \log N \right)
\cap \left( \max _{1 \leq i \leq T(N)} X_i < \sqrt{N}\right)
\right\} +o(1)\]
\[ \leq 
{\bf P} \left\{\sqrt{N} \, a \log N > bN 
\right\} +o(1)=o(1).\]
Substituting this bound into (\ref{C1}) we immediately get (\ref{C}). \hfill$\Box$

\bigskip

Conditionally on ${\mathcal{C}}_1 =V(\widetilde{\mathcal{C}}_1)$ we have
\begin{equation}\label{A26}
\begin{array}{ll}
\frac{C_1}{N} & 
= \frac{1}{N} \sum_{i=1}^{T(N)} X_i {\bf 1}\{
v_i \in \widetilde{\mathcal{C}}_1 \} \\
\\
& = \frac{1}{N} \sum_{i=1}^{T(N)}\sum_{k=1}^{N} k 
{\bf 1}
\{
X_i=k \}  {{\bf 1}}\{
v_i \in \widetilde{\mathcal{C}}_1 \} \\
\\
& = \frac{T(N)}{N} \sum_{k=1}^{N} k \frac{1}{T(N)}
\#\{
v_i \in \widetilde{\mathcal{C}}_1 : X_i=k\}. 
\end{array}
\end{equation}
Note that Theorem 9.10 from \cite{BJR} (together with Proposition \ref{P1} in our case)
 implies  that 
\begin{equation}\label{A32}
\nu_N(k):=\frac{1}{T(N)} \# \{
v_i \in \widetilde{\mathcal{C}}_1 (N): X_i=k\} \stackrel{P}{\rightarrow}
 {
  \rho}(\kappa ; k) \mu (k)
 \end{equation}
for each $k \geq 1$.

We shall prove below that also 
\begin{equation}\label{A27}
W_N:=\sum_{k=1}^{N} k\nu_N(k)\stackrel{P}{\rightarrow}
 \sum_{k=1}^{\infty} k {
  \rho}(\kappa ; k) \mu (k) =:  \beta \, {\bf E}X.
\end{equation}
Observe that ${
  \rho}(\kappa ; k)$ is the maximal solution to (\ref{rho}), therefore
$\beta $ is the maximal solution to
\[\beta = \frac{1}{{\bf E}X}\sum_{k=1}^{\infty} k {
  \rho}(\kappa ; k) \mu (k)=\frac{1}{{\bf E}X}\sum_{k=1}^{\infty} k 
\left(
1- e^{-\sum_{y=1}^{\infty} \kappa
 (k,y) \mu(y)\rho(\kappa;y) }\right)
\mu (k)\]
\[= 1- \frac{1}{{\bf E}X}{\bf E}\Big( Xe^{-cX\beta} \Big).
\]
 This proves that 
$\beta$  is the maximal root of (\ref{be}).
Then (\ref{A27}) together with Proposition \ref{P1},
which states that 
$T(N)/N \stackrel{a.s.}{\rightarrow} 1/{\bf E}X$, will allow us to
derive from 
(\ref{A26}) that for any positive $\varepsilon $
\[
{\bf P}
 \Big\{ |\frac{ C_{1}\Big(G_N(c,q)\Big)}{N} - \beta | > \varepsilon \mid 
{\mathcal{C}}_1 =V(\widetilde{\mathcal{C}}_1) \Big\}
\rightarrow 0
\]
as $N \rightarrow \infty$. This combined with Lemma \ref{P4} would
immediately imply
\begin{equation}\label{C1/N}
\frac{ C_{1}\Big(G_N(c,q)\Big)}{N} \stackrel{P}{\rightarrow}  \beta ,
\end{equation}
and hence the statement of the theorem follows.

Now we are left with proving (\ref{A27}). 
For any $1 \leq R < N$ write $W_N:=W_N^{R}+w_N^{R}$, where
\[W_N^{R}:=\sum_{k=1}^{R} k\nu_N(k), \ \ 
w_N^{R}:=\sum_{k=R+1}^{N} k\nu_N(k).
\]
By (\ref{A32}) we have for any fixed $R\geq 1$
\begin{equation}\label{A33}
W_N^{R} \ \stackrel{P}{\rightarrow} \ 
 \sum_{k=1}^{R} k {
  \rho}(\kappa ; k) \mu (k) 
\end{equation}
 as $N\rightarrow \infty$. Consider $w_N^{R}$. Note 
that for any $k\geq 1$
\begin{equation}\label{A34}
{\bf E}\nu _N(k)={\bf E}{\bf E}\Big\{\nu _N(k) \mid T(N)\Big\} \leq
{\bf E}
\frac{1}{T(N)} \sum_{i=1}^{T(N)} {\bf P}\{
 X_i=k\mid T(N)\}.  
\end{equation}
Using events ${\cal A}_{\delta, N}$ with
bound (\ref{AA})  and
Proposition \ref{P1} 
we obtain from (\ref{A34}) for any fixed $0<\delta <1/(2{\bf E} X)$
\[
{\bf E}\nu _N(k) \leq
{\bf E}
\frac{1}{T(N)} \sum_{i=1}^{T(N)} {\bf P}\{
 X_i=k\mid T(N)\} {\bf 1}\{{\cal A}_{\delta, N}\} +{\bf P}\{\overline{{\cal A}_{\delta, N}}\}  
\]
\[ \leq  \frac{(1+\delta {\bf E} X) }{(1-\delta {\bf E} X)} \
{\bf P}\{
 X_1=k \} (1+o(1))+{\bf P}\{\overline{{\cal A}_{\delta, N}}\} . 
\]
Bound (\ref{AA}) allows us to derive from here that
\begin{equation}\label{A35}
{\bf E}\nu _N(k) \leq A_1 (\mu(k)+e^{-a_1 N})
\end{equation}
for some positive constants $A_1$ and $a_1$ independent of $k$ and $N$.
This yields
\begin{equation}\label{A36}
{\bf E} w_N^{R}=\sum_{k=R+1}^{N} k 
{\bf E}\nu_N(k) \leq A_2 e^{-a_2 R}
\end{equation}
for some positive constants $A_2$ and $a_2$. 

Now for any $\varepsilon >0$ we can choose $R$ so that 
\[\sum_{k=R+1}^{\infty} k {
  \rho}(\kappa ; k) \mu (k) < \varepsilon /3,\]
and then we have
\begin{equation}\label{A37}
{\bf P} \{| W_N-\sum_{k=1}^{\infty} k {
  \rho}(\kappa ; k) \mu (k) |>\varepsilon\}
\end{equation}
\[={\bf P} \{|( W_N^R -\sum_{k=1}^{R} k {
  \rho}(\kappa ; k) \mu (k))
 +w_N^R-
\sum_{k=R+1}^{\infty} k {
  \rho}(\kappa ; k) \mu (k) |>\varepsilon\}\]
\[\leq {\bf P} \{| W_N^R -\sum_{k=1}^{R} k {
  \rho}(\kappa ; k) \mu (k)|>\varepsilon /3\}
 +{\bf P} \{w_N^R>\varepsilon /3\}. \]
Markov's inequality together with bound (\ref{A36}) gives us 
\begin{equation}\label{A38}
{\bf P} \{w_N^R>\varepsilon /3\} \leq \frac{3 {\bf E} w_N^{R}}{\varepsilon}
\leq \frac{3 A_2 e^{-a_2 R}}{\varepsilon}.
\end{equation}
Making use of (\ref{A38}) and (\ref{A33}) we immediately derive from 
(\ref{A37}) 
\begin{equation}\label{A39}
{\bf P} \{| W_N-\sum_{k=1}^{\infty} k {
  \rho}(\kappa ; k) \mu (k) |>\varepsilon\}
\leq o(1) + \frac{3 A_2 e^{-a_2 R}}{\varepsilon}
\end{equation}
as $N\rightarrow \infty$. Hence, for any given positive $\varepsilon$ and
$\varepsilon _0$ we can choose finite $R$ so large that  
\begin{equation}\label{A40}
\lim _{N \rightarrow \infty}{\bf P} \{| W_N-\sum_{k=1}^{\infty} k {
  \rho}(\kappa ; k) \mu (k) |>\varepsilon\}
< \varepsilon_0.
\end{equation}
This clearly proves statement (\ref{A27}), and therefore finishes the proof of 
the theorem.
\hfill$\Box$

\end{document}